\documentclass [twoside,10pt]{article}
\usepackage{graphicx}
\usepackage{amssymb}
\usepackage{amsthm}
\newtheorem{lem}{Lemma}

\newcommand{\diam}{{\rm diam}}

\newcommand{\sh}{{\rm sinh \,}}
\newcommand{\ch}{{\rm cosh \,}}

\newcommand{\bd}{{\rm bd}}
\newcommand{\conv}{{\rm conv}}

\date{}

\begin{document}

\title{$\!\!\!\!\!\!\!\!\!\!\!\!\!$Spherical quadrilateral with three right angles and its application for diameter of extreme points of a convex body}

\maketitle

\thispagestyle{empty}

\vskip-1cm

\noindent
{\author{Marek Lassak}}

\pagestyle{myheadings} \markboth{\centerline {Marek Lassak}}{\centerline {Spherical quadrilateral with three right angles and its application}}

\maketitle
\vskip 0.5cm

\noindent
{\bf Abstract}.
We prove a theorem on the relationships between the lengths of sides of a spherical quadrilateral with three right angles.
They are analogous to the relationships in the Lambert quadrilateral in the hyperbolic plane.
We apply this theorem in the proof of our second theorem that if $C$ is a two-dimensional spherical convex body of diameter $\delta \in (\frac{1}{2}\pi,\pi)$, then the diameter of the set of extreme points of $C$ is at least $2 \arccos \big(\frac{1}{4}(\cos \delta + \sqrt {\cos^2 \delta +8})\big)$.
This estimate cannot be improved.

\baselineskip 16.8pt



\medskip

\date{} 

\maketitle

\section{Introduction}

Denote by $S^n$ be the unit sphere of the $(n+1)$-dimensional Euclidean space $E^{n+1}$. 
In this article we mostly deal with $S^2$.

Let us recall some needed notions from the spherical geometry.
By a {\it great circle} of $S^2$ we mean the intersection of $S^2$ with any two-dimensional linear subspace of $E^{3}$. 
By a pair of {\it antipodes} of $S^2$ we understand any pair of points being the intersection of $S^2$ with a one-dimensional subspace of $E^3$.
Observe that if two different points $p, q \in S^2$ are not antipodes, there is exactly one great circle passing through them.
By the {\it arc} $pq$ connecting these points we mean the shorter part of the great circle through $p$ and $q$, and 
by the {\it distance} $|pq|$ of $p$ and $q$ we understand the length of the arc $pq$. 

A subset of $S^2$ is called {\it convex} if it does not contain any pair of antipodes of $S^2$ and if together with every two points it contains the arc connecting them.
A closed convex set $C \subset S^2$ with non-empty interior is called a {\it convex body}.
Its boundary is denoted by $\bd (C)$.
A point of $C$ is called {\it extreme} if it does not belong to the relative interior of any arc with end-points in $C$.
The set of extreme points of $C$ is denoted by $E(C)$.
Let us add that convexity in $S^n$ is considered in very many papers and monographs (see [1--12]).

The set of points of $S^2$ in the distance at most $\rho$, where $0 < \rho \leq \frac{\pi}{2}$, from a fixed point of $S^2$ is called a {\it spherical disk} of {\it radius} $\rho$.  
Spherical disks of radius $\frac{\pi}{2}$ are called {\it hemispheres}.
Two hemispheres whose centers are antipodes are called {\it opposite hemispheres}.
The set of points of a great circle of $S^2$ which are at a distance at most $\frac{\pi}{2}$ from a fixed point $p$ of this great circle is called a {\it semicircle}.
We call $p$ the {\it center} of this semicircle.

If non-opposite hemispheres $G$ and $H$ are different, then $L = G \cap H$ is called a {\it lune}. 
As in \cite{L1} and \cite{L2}, the semicircles bounding $L$ and contained in $G$ and $H$, respectively, are denoted by $G/H$ and $H/G$. 
The {\it thickness} $\Delta (L)$ of $L$ is defined as the distance between the centers of $G/H$ and $H/G$.  
The two common points of the semicircles bounding $L$ are called {\it corners} of~$L$. 

Recall a lemma from  \cite{LasMus-FAS}: 
{\it for every convex body $C \subset S^2$ of diameter at most $\frac{\pi}{2}$ we have ${\rm diam} (E(C)) = {\rm diam}(C)$.} 
In the present paper we consider an analogous property for every convex body $C \subset S^2$ of diameter over $\frac{\pi}{2}$.
Namely, it is proved in Theorem 2 that $ \diam(E(C)) \geq 2 \arccos \big(\frac{1}{4}(\cos \delta + \sqrt {\cos^2 \! \delta +8})\big)$ for every convex body $C \subset S^2$ of diameter over $\frac{\pi}{2}$.
In order to overcome some troubles in its proof, it is necessary to find a length of a side of a quadrilateral with three right angles.

We perform the above task in Theorem 1, where we present some relationships between the lengths of sides of a spherical quadrilateral with three right angles, similar to the relationships in the Lambert quadrilateral in the hyperbolic plane.

\section{Some properties of lunes} 

The following function, where $\delta \in (\frac{1}{2}\pi, \pi)$, plays an important role in our research:

\vskip-0.5cm
$$\varphi(\delta) = \arccos \big(\frac{1}{4}(\cos \delta + \sqrt {\cos^2 \! \delta +8})\big).$$

\begin{lem} \label{one-valued}
{\it The function $\varphi(\delta)$ is increasing and thus one-valued in the interval $(\frac{1}{2}\pi, \pi)$.
The set of its values is the interval $(\frac{1}{4} \pi, \frac{1}{3} \pi)$.
We have $\cos \delta = \frac{2\cos^2 \! \varphi - 1}{\cos \varphi}$ for $\varphi$ from this interval, where 
$\varphi(\delta)$ is shortly denoted by~$\varphi$.}
\end{lem}

\begin{proof}
The derivative of $\; \cos \delta + \sqrt {\cos^2 \! \delta +8}\;$ is $\;-\sin \delta \;\!(1 + \cos \delta /\! \sqrt{\cos^2 \delta +8})$.
Since $\sqrt {\cos^2 \! \delta +8} > -\cos \delta \;$ for $\; \delta \in (\frac{1}{2}\pi, \pi)$, we conclude that this derivative is negative in the interval $(\frac{1}{2}\pi, \pi)$.
Hence $\cos \delta + \sqrt {\cos^2 \! \delta +8}$ is a decreasing function of $\delta$ in the interval $(\frac{1}{2}\pi, \pi)$.
So $\varphi(\delta)$ is increasing there.
Hence it is one-valued.  
Moreover, from $\varphi (\frac{1}{2}\pi) = \frac{1}{4} \pi$ and $\varphi (\pi) = \frac{1}{3} \pi$, by the continuity of $\varphi (\delta)$ we get the second thesis.
The third one results by a simple calculation.
\end{proof}

Here is a lemma analogous to Lemma from \cite{LasMus-FAS}, where a lune of thickness at most $\frac{\pi}{2}$ was considered.
We keep the same notation.

\begin{lem} \label{min}
Let $L \subset S^2$ be a lune of thickness over $\frac{\pi}{2}$ whose bounding semicircles are $Q$ and $Q'$.
For every $u, v, z$ in $Q$ such that $v \in uz$ and for every $q \in Q'$ we have $|qv| \geq \min \{|qu|, |qz| \}$.
\end{lem}

\begin{proof}
Clearly, the farthest point $p \in Q$ to $q$ is unique.
Observe that for $x\in Q$ the distance $|qx|$ decreases as the distance $|px|$ increases. 
This easily implies the assertion of our lemma. 
\end{proof}

\vskip0.1cm
\noindent
{\bf Proposition.} 
{\it Denote by $L$ a lune of thickness $\delta \in (\frac{\pi}{2}, \pi)$. and by $G$ and $H$ the semicircles bounding it.
Denote by $g$ the center of $G$ and by $h$ the center of $H$.
Let $r$ and $s$ be the two points of $G$ in the distance $\varphi(\delta) = \arccos \big(\frac{1}{4} (\cos \delta + \sqrt {\cos^2 \! \delta +8})\big)$ from $g$. 
For any $t \in rs$, by $t'$ we mean the intersection of $H$ with the great circle through $t$ orthogonal to $H$.

I. The distances $|rh|$ and $|sh|$ are equal $2\varphi(\delta)$.
Moreover, $|th| \geq 2\varphi(\delta)$ for every $t \in rs$.

II. The triangle $rsh$ is equilateral with the length $2\varphi(\delta)$ of its sides.

III. If $l \in ht'$, then $|tl| \geq 2\varphi(\delta)$.}

\begin{proof}
Denote by $c_r$ this corner of $L$ which is closer to $r$, and by $c_s$ this one closer to $s$.
Of course, $|gh| = \delta$ and $|hc_r| = |hc_s| = \frac{1}{2}\pi$.
Observe that $|hx|$ grows as $x$ moves from $c_r$ (or from $c_s$) to $g$ on $G$. 
Hence there exists a unique point $y$ on $gc_r$ such that $|yh| = 2|yg|$.
Also on $gc_s$, there is exactly one.
Clearly, $hgy$ is a right triangle with the right angle at $g$.
Of course, $\cos |yh| = \cos |gy|  \cos |gh|$. 
Denote $|gy|$ by $\psi$.
Then $|yh| =2\psi$.
Hence $\cos 2\psi = \cos \psi  \cos \delta$.
From $|gh| = \delta$ and the trigonometric formula $\cos 2\psi = 2\cos^2 \! \psi -1$ we obtain the equation $2\cos^2 \! \psi -1 = \cos \psi  \cos \delta$ with the unknown $\psi$. 
Its solutions is $\psi = \arccos \big(\frac{1}{4} (\cos \delta + \sqrt {\cos^2 \! \delta +8})\big)$.
We see that the required point $y$ can be only at $r$ or $s$.
This confirms the first thesis of Part I.
The second thesis of Part I results from the earlier observed property that $|hx|$ grows as $x$ moves from $c_r$ (or from $c_s$) to $g$ on $G$. 

From Part I we obtain Part II.

Let us show Part III.
Observe that $|tt'|$ is at least $|gh| = \delta$, which implies that it is over $\frac{1}{2}\pi$.
Clearly, $|tl|$ grows as $l$ moves from $h$ to $t'$ on $ht'$.
So $|tl| > \frac{1}{2}\pi$ for every $l \in ht'$.
Hence $|tl| > |th|$.
Consequently, by the second statement of Part I we get $|tl| \geq 2 \varphi(\delta)$.
\end{proof}

By the way, having in mind Part II of Proposition, the author expects that a larger equilateral triangle cannot be a subset of $L$.

By the second statement of Part I of Proposition we conclude the following claim.

\vskip0.2cm
\noindent
{\bf Claim.}
{\it For every $\delta \in (\frac{\pi}{2}, \pi)$ we have $\delta \geq 2\varphi(\delta)$.}

\section{Quadrilaterals with three right angles}

Let $Q= abcd$ be a spherical quadrilateral with right angles at $a, b$ and $c$.
Put $\kappa = |ab|$, $\lambda = |bc|$, $\mu =|cd|$ and $\nu = |da|$.
Recall that for the alike hyperbolic Lambert quadrilateral $abcd$ we have
we have $\sh \mu = \sh \kappa \; \ch \nu$, $\tanh \mu = \ch \lambda \tanh \kappa$, $\sh \nu = \sh \lambda \; \ch\mu$, $\tanh \nu = \ch \kappa \tanh \lambda$.
Here is a theorem in which we present some relationships in the quadrilateral $Q$ analogous to the ones in the hyperbolic Lambert quadrilateral.

\newpage
\noindent
{\bf Theorem 1.} 
{\it In the spherical quadrilateral $Q$ we have}

(1) $\sin \mu =  \sin \kappa \cos \nu,$ 

(2) ${\tan \mu} =  \tan \kappa \cos \lambda,$ 

(3) $\cos\nu = \sqrt{\cos^2 \! \mu \cos^2 \! \lambda + \sin^2 \! \mu},$ 

(4) $\cos \nu = \frac{\cos \lambda}{\sqrt{1 - \sin^2 \! \lambda \sin^2 \! \kappa}}.$

\begin{proof}
Denote $|bd|$ by $\xi$ (see Figure 1). 
From the right triangle $bcd$ we get $\cos \xi = \cos \mu  \cos \lambda$.
From the right triangle $dab$ we get $\cos \xi = \cos \nu \cos \kappa  $.
Hence $\cos \mu  \cos \lambda =  \cos \nu \cos \kappa $.

\begin{figure}[htbp]   
\hskip2.1cm     
\includegraphics[width=10.16cm, height=6.8cm]{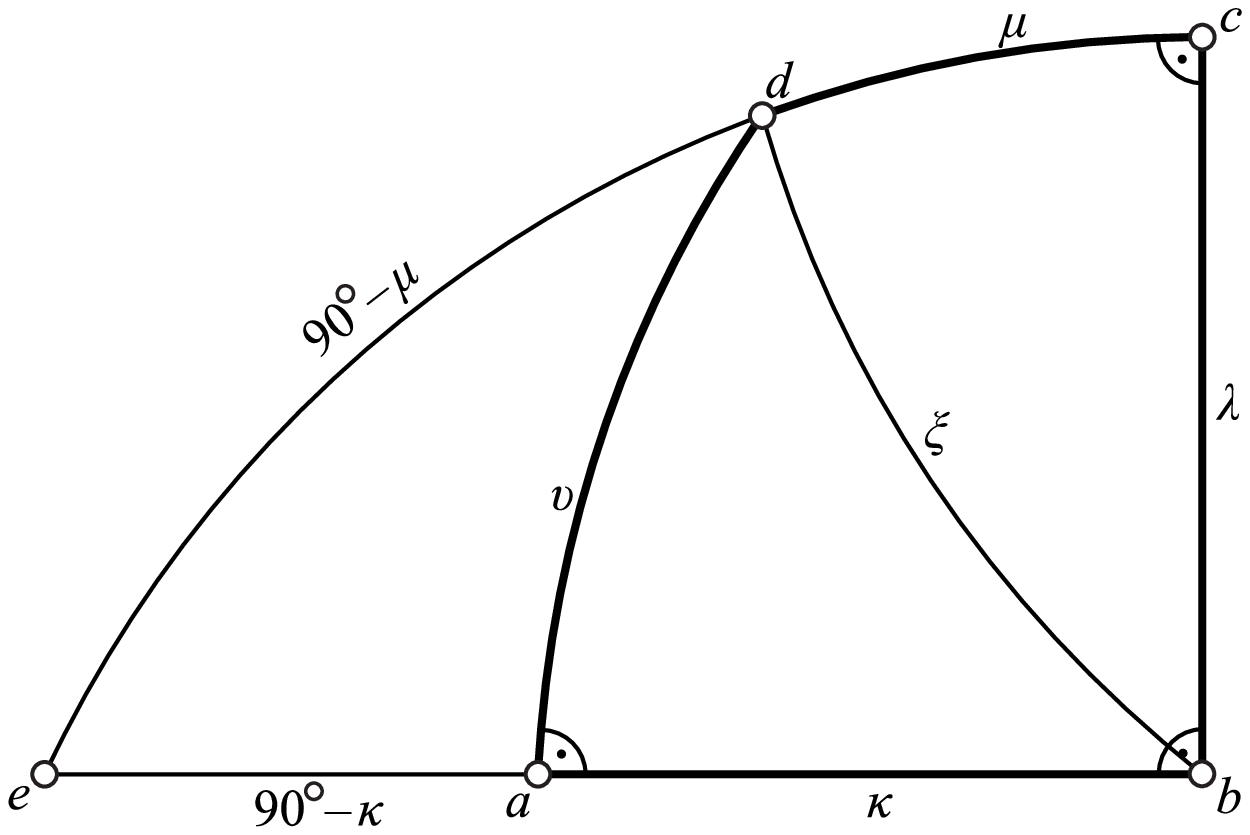}

\vskip0.3cm 
\centerline
{{\bf Fig. 1.} Spherical quadrilateral with three right angles} 
\vskip-0.1cm 
\end{figure}

Prolong the arcs $cd$ and $ba$ up to the intersection point $e$ (we mean that $d \in ce$ and $a \in eb$).
From the right triangle $ead$ we see that $\cos |de| = \cos |ea| \cos |ad|$.
Moreover, since $ce$ and $ba$ are perpendicular to the great circle containing $bc$, we have $|ce| = 90^\circ$ and $|be| = 90^\circ$.
Consequently, from $\cos |de| = \cos (90^\circ - \mu) = \sin \mu$ and $\cos |ea| = \cos (90^\circ - \kappa) = \sin \kappa$ we obtain~(1).

The preceding two paragraphs imply (2).
Still it is sufficient to divide  left and right sides of $\sin \mu = \cos \nu \sin \kappa$ by left and right sides of $\cos \mu  \cos \lambda = \cos \nu \cos \kappa$, respectively, which gives $\frac{\tan \mu}{\cos \lambda} = \tan \kappa$.
This is equivalent to (2).

Let us show (3).
From (2) we get $\tan^2 \! \mu = \tan^2 \! \kappa \cos^2 \! \lambda$
which equals $\frac{\sin^2 \!\kappa}{1-\sin^2 \!\kappa}\cos^2\lambda$.
Applying (1) we substitute here $\sin \kappa$ by  $\frac{\sin\mu}{\cos \nu}$.
We get $\tan^2 \!\mu = \frac{\sin^2 \!\mu}{\cos^2 \! \nu - \sin^2 \! \mu} \cos^2 \!\lambda$.
Consequently, $\cos^2 \! \nu = \cos^2 \!\mu\cos^2 \!\lambda +\sin^2 \!\mu$, which gives (3).

Finally, we intend to show (4).
We start from the squared (3).
Since $\cos^2 \!\mu = 1-\sin^2 \!\mu$, we have
$\cos^2 \!\nu = (1-\sin^2 \!\mu)\cos^2 \!\lambda +\sin^2 \!\mu = \sin^2 \!\mu (1-\cos^2 \!\lambda) +\cos^2 \!\lambda = \sin^2 \!\mu\sin^2 \!\lambda +\cos^2 \!\lambda$.
Having in mind (1), let us substitute here $\sin^2 \!\mu = \sin^2 \!\kappa\cos^2 \!\nu$.
We get $\cos^2 \!\nu = \sin^2 \!\kappa \cos^2 \!\nu \sin^2 \! \lambda + \cos^2 \! \lambda$. 
Hence $\cos^2 \!\nu (1-\sin^2 \!\kappa\sin^2 \!\lambda) = \cos^2 \!\lambda$, which implies (4).
\end{proof}

\section{Comparison between the diameters of $C$ and $E(C)$}

\vskip0.1cm
\noindent
{\bf Theorem 2.} {\it For every convex body $C \subset S^2$ of diameter $\delta \in (\frac{\pi}{2}, \pi)$ we have} 
\vskip-0.25cm
$$ \diam(E(C)) \geq 2 \arccos \big(\frac{1}{4}(\cos \delta + \sqrt {\cos^2 \! \delta +8})\big).$$

\noindent
{\it This inequality cannot be improved for every $\delta$.}

\vskip-0.25cm
\begin{proof} 
Shortly saying, the inequality in our theorem is $ \diam(E(C)) \geq 2\varphi(\delta)$.

By compactness arguments, the diameter $\delta$ of $C$ is realized for a pair $g, j$ of boundary points of $C$.
Clearly, $g, j \in \bd(C)$.

By Lemma 7 of \cite{LasMus-AEQ} the semicircle $T$ of radius $\delta$ centered at $g$ and orthogonal to $gj$ at $g$ supports $C$. 
Analogously, the semicircle $W$ $\delta$ centered at $j$ and orthogonal to $gj$ at $j$ supports $C$. 
Denote by $o$ the center of the lune $L = \conv (T \cup W)$, this is the midpoint of $gj$.
Clearly, $C \subset L$ and $L$ has thickness $|gj| = \delta$. 
Consider three cases.

\vskip0.25cm
Case 1, when $g, j \in E(C).$ 

\vskip0.1cm
Since $|gj| = \delta$, by Claim we get the thesis of our theorem.

\vskip0.2cm
Case 2, when $j \in E(C)$ and $g \not \in E(C)$ (or analogously, when $j \not \in E(C)$ and $g \in E(C)$).

\vskip0.1cm
Since $g \not \in E(C)$, there are $f, h \in E(C)$ different from $g$ such that $g \in fh$. 
Clearly, these two points belong to the semicircle $T$.
We intend to show that $\max \{|jf|, |jh|, |fh| \} \geq 2\varphi(\delta)$.
This will imply $\diam(E(C)) \geq 2\varphi(\delta)$, what is our aim. 

Observe that for every position of $f, h$ in $T$ at least one of the following three sub-cases holds true.
In each we obtain a consequence.

\vskip0.1cm
(2.1) Assume that $|gf| \leq \varphi(\delta)$. 
Then $|jf| \geq 2\varphi(\delta)$, as it follows by the second thesis of Part I of Proposition.

(2.2) Assume that $|gh| \leq \varphi(\delta)$. 
Then $|jh| \geq 2\varphi(\delta)$, as it follows by the second thesis of Part I of Proposition.

(2.3) Assume that $|gf| \geq \varphi(\delta)$ and  $|gh| \geq \varphi(\delta)$.
Then $|fh| \geq 2\varphi(\delta)$.

\vskip0.15cm
As explained earlier in Case 2, from the obtained consequences we conclude that $\diam(E(C)) \geq 2\varphi(\delta)$.

\vskip0.2cm
Case 3, when $g \not \in E(C)$ and $j \not \in E(C)$.

\vskip0.1cm
As in the preceding case, there are points $f, h \in E(C)$ different from $g$ such that $g \in fh$.
Clearly, again $f, h \in W$.
We do not lose the generality assuming that $f, g, h$ lie counterclockwise on $T$.
Since $j \not \in E(C)$, there are $i, k \in E(C)$ different from $j$ such that $j \in ik$ and $i, j, k$ lie counterclockwise on $W$.
Our aim is to show that $\diam \{f, h, i, k\} \geq 2 \varphi(\delta)$.
Still it implies  $\diam(E(C)) \geq 2 \varphi(\delta)$.

\vskip0.1cm
For every $f, h, i, k$ at least one of the following four conditions is true.

(3.1) \ $|ji| \geq |gf|$ and $|jk| \geq |gh|$,

(3.2) \ $|ji| \leq |gf|$ and $|jk| \leq |gh|$,

(3.3) \ $|ji| \leq |gf|$ and $|jk| \geq |gh|$,

(3.4) \ $|ji| \geq |gf|$ and $|jk| \leq |gh|$.

\vskip0.1cm 
In (3.1) consider the following three disjoint possibilities.
In each we present a consequence.

(3.1.a) Let $|jk| \leq \varphi(\delta)$.
Denote by $g^*$ the center of the great circle containing $T$ and by $k'$ the projection of $k$ onto $T$. 
Clearly, $g^* \in kk'$.
Let us apply Part III of Proposition taking $k, k', h$ in the parts of $t, t',l$, respectively.
We conclude that $h \in gk'$ and get $|hk| \geq 2 \varphi(\delta)$.

\vskip0.1cm
(3.1.b) Let $|ji| \leq \varphi(\delta)$.
Analogously as in (3.1) we get $|fi| \geq 2\varphi(\delta)$.
  
\vskip0.1cm
(3.1.c) 
Let $|ji| > \varphi(\delta)$ and $|jk| > \varphi(\delta)$.
Then $|ik|  > 2\varphi(\delta)$.

\vskip0.1cm
We see that always $\max \{|hk|, |fi|, |ik| \} \geq 2\varphi(\delta)$.
Hence $\diam(E(C)) \geq  2\varphi(\delta)$. 

\vskip0.1cm
If (3.2) is true, we deal analogously as in (3.1).

\vskip0.3cm
Consider (3.3).

If $|gf| \leq \varphi(\delta)$, then also $|ji| \leq \varphi(\delta)$.
Taking into account the projection $f'$ of $f$ onto $W$ and applying Part III of Proposition we obtain $|fi| \geq 2\varphi(\delta)$. 

Analogously, if $|jk| \leq \varphi(\delta)$, then also $|gh| \leq \varphi(\delta)$ and by Part III of Proposition we obtain 
$|hk| \geq 2\varphi(\delta)$.

If both the above assumptions (in the two preceding paragraphs) are not fulfilled , so if  $|gf| > \varphi(\delta)$ and $|jk| > \varphi(\delta)$, then we consider the following three possibilities. 

(3.3.1) when $|gh| > \varphi(\delta)$.
Clearly $|fh| > 2\varphi(\delta)$. 

(3.3.2) when $|ji| > \varphi(\delta)$.
Clearly $|fh| > 2\varphi(\delta)$. 

(3.3.3) when $|gh| \leq \varphi(\delta)$ and $|ji| \leq \varphi(\delta)$.  
Below we consider three (not disjoint) possibilities.

\vskip0.1cm
(3.3.3.a) when $\frac{1}{2}\varphi(\delta) \leq |gh| \leq \varphi(\delta)$.
By Lemma 1 we have $\varphi(\delta) < \frac {1}{3}\pi$ provided $\frac{1}{2}\pi <  \delta < \pi$. 
Hence $\frac{1}{2}\pi > \frac{3}{2}\varphi(\delta)$.
Thus 

$$\frac{1}{2}\varphi(\delta) + \frac{1}{2}\pi > 2\varphi(\delta). \eqno (*)$$

Moreover, let us show that 

\vskip-0.5cm
$$|hk| >  \frac{1}{2}\varphi(\delta) + \frac{1}{2}\pi. \eqno (**)$$

\vskip0.1cm
Denote by $c_f$ the corner of $L$ closer to $f$ and by $c_h$ the corner of $L$ closer to $h$.
Provide the orthogonal great circle to $W$ through $h$.
Clearly, it passes through $jc_f$.
Denote by $h'$ the point of intersection with $jc_f$. 
Since $k \in jc_f$, consider the following two possibilities. 

Look at the first possibility when $k \in h'c_f$. 
In particular, if $k$ is at $c_f$, then $|hc_f| \geq \frac{1}{2}\varphi(\delta) + \frac{1}{2}\pi$ is true since $hc_f$ is a part of $T$ and $|gh| \geq \frac{1}{2}\varphi(\delta)$. 
Moreover, since $|hh'| > \frac{1}{2}\pi$, then also $|hx| > \frac{1}{2}\pi$ for every $x \in c_fh'$, which implies that if $x$ moves on $c_fh'$ from $c_f$ to $h'$, then $|hx|$ grows and consequently $|hk| \geq \frac{1}{2}\varphi(\delta) + \frac{1}{2}\pi$. 
Hence if $k \in h'c_f$ then also $|hk| \geq \frac{1}{2}\varphi(\delta) + \frac{1}{2}\pi$. 
Thus $|hk| > \frac{1}{2}\varphi(\delta) + \frac{1}{2}\pi$, i.e., ($**$) is true.  

Consider now the second possibility when $k \in jh'$. 
Observe that since $|hh'| > \frac{1}{2}\pi$, then also $|hk| > \frac{1}{2}\pi$. 
Thus $|hx|$ grows as $x$ moves from $j$ to $h'$ on $jh'$. 
Hence $|hk| > |hj|$.
Moreover, by the assumption $|gh| < \varphi(\delta)$ of (3.3.3) and Part III of Proposition we get $|hk| \geq 2 \varphi(\delta)$.
So $|hk| > 2 \varphi(\delta)$.

We see that ($**$) holds true for both the possibilities.

From ($*$) and ($**$) we obtain $|hk| > 2\varphi(\delta)$.

\vskip0.1cm
(3.3.3.b) when $\frac{1}{2}\varphi(\delta) \leq |ji| \leq \varphi(\delta)$.
We deal analogously as in (3.3.3.a). 

\vskip0.1cm
(3.3.3.c) when none of the inequalities $\frac{1}{2}\varphi(\delta) \leq |gh| \leq \varphi(\delta)$ and $\frac{1}{2}\varphi(\delta) \leq |ji| \leq \varphi(\delta)$ is true.
We have $|gh| < \frac{1}{2}\varphi(\delta)$ and $|ji| < \frac{1}{2}\varphi(\delta)$.

We may assume that, for instance,  $|ji| \leq |gh|$ (see Figure 2).

Take into account two possibilities:

(3.3.3.c') when $|gh| \geq 2\varphi(\delta) - \frac{\pi}{2}$, 

(3.3.3.c'') when $|gh| < 2\varphi(\delta) - \frac{\pi}{2}$.

\vskip0.1cm
Consider (3.3.3.c').

By Lemma \ref{min} we have $|hk| \geq \min \{|hc_f|, |hj|\} \geq 2\varphi(\delta)$. 
The reason is that $|hc_f| = |hg| + |gc_f| \geq [2\varphi(\delta) - \frac{1}{2}\varphi(\delta)] + \frac{1}{2}\varphi(\delta) = 2\varphi(\delta)$ and that $|hj| \geq |sj| = 2\varphi(\delta)$, where $s \in gc_h$ is the point in the distance $\varphi(\delta)$ from~$g$ (compare Proposition and especially its Part II).

\vskip0.1cm
Consider (3.3.3.c'').

By $h^+$ denote the point being the intersection of $W$ with the great circle through $h$ orthogonal to $c_fc_h$ (so $h^+$ is the ``symmetric" point to $h$ with respect to $c_fc_h$).

We intend to show that $|hi| \geq 2\varphi(\delta)$.
Since $|hh^+| \leq |hi|$, it is sufficient to show that $|hh^+| \geq 2\varphi(\delta)$.
Thus it is sufficient to show that $|hm| \geq \varphi(\delta)$, where $m$ denotes the midpoint of $hh+$.

First let us show that $|hk| \geq 2\varphi(\delta)$ in the particular worst situation when $|gh| = 2\varphi(\delta) - \frac{\pi}{2}$.
Since the quadrangle $hmog$ has right angles at $m, o$ and $g$, let us apply (3) from Theorem 1, where $m=a, o=b, g=c$ and $h=d$. 
So in the terms of Theorem 1 we assume that $\mu = 2\varphi - \frac{\pi}{2}$
and we intend to show that $\nu \geq \varphi(\delta)$.
In other words, we intend to show that $\cos \nu \leq \cos \varphi(\delta)$.
Having in mind (3) from Theorem 1, we intend to show that

\vskip-0.1cm
$$\sqrt{\cos^2 (2\varphi(\delta) -\frac{\pi}{2}) \cos^2 \! \lambda + \sin^2 \! (2\varphi(\delta)- \frac{\pi}{2})} \leq \cos \varphi(\delta). \eqno (***)$$

\begin{figure}[htbp]        
\hskip2.2cm
\includegraphics[width=10.13cm, height=7.55cm]{lb241202}
\vskip0.2cm 
\centerline
{{\bf Fig. 2.} Illustration to the proof of (3.3.3.c) in Theorem 2} 
\vskip-0.01cm 
\end{figure}

\vskip-0.2cm
In other words, we intend to show that
$\sqrt{\sin^2 (2\varphi(\delta)) \cos^2 \! \lambda + \cos^2 (2\varphi(\delta))} \leq \cos \varphi(\delta).$

We deal here with an arbitrary $\delta$ such that $\frac{1}{2}\pi < \delta < \pi$, so with $\lambda \in (\frac{1}{4}\pi, \frac{1}{2}\pi)$ and $\varphi(\delta) \in (\frac{1}{4}\pi, \frac{1}{3}\pi)$, what Lemma \ref{one-valued} says.
Since $\lambda = \frac{1}{2}\delta$, we are planning to show that

$$\sqrt{\sin^2 2\varphi \cos^2 \! \frac{\delta}{2} +  \cos^2 \! 2\varphi} \leq \cos \varphi,$$ 

\noindent
where for simplicity we write $\varphi$ in place of $\varphi(\delta)$. 
This is equivalent to
$\sin^2 2\varphi \cos^2 \! \frac{\delta}{2} + \cos^2 \! 2\varphi \leq \cos^2 \varphi$. 

This is equivalent to 
$\sin^2 2\varphi \cdot \frac{1+\cos\delta}{2} +\cos^2 2\varphi \leq \cos^2\varphi$ 
and by Lemma \ref{one-valued} also to 
$(1-\cos^2 2\varphi) \cdot \frac{1+ \frac {2\cos^2 \varphi -1}{\cos \varphi}}{2} + \cos^2 2\varphi \leq \cos^2\varphi$.
Having in mind that $\cos 2\varphi = 2\cos^2 \varphi -1$ we see that this is equivalent to
$\big(1- (2\cos^2\varphi -1)^2 \big) \cdot \frac{1+ \frac {2\cos^2 \varphi -1}{\cos \varphi}}{2} + (2\cos^2 \varphi -1)^2 \leq \cos^2\varphi$.

Putting $x = \cos \varphi$, we get $(1- (2x^2-1)^2) \cdot (1 + \frac{2x^2 - 1}{x})/2 + (2x^2-1)^2  \leq x^2$, i.e., that $(2x-1)(x+1)(x-1)(1-2x^2) \leq 0$.
In particular, this inequality is true for $x \in [\frac{1}{2}, \frac{1}{2}\sqrt 2]$. 
In other words, for $x \in [\cos \frac{\pi}{3}, \cos \frac{\pi}{4}]$.
Hence (***) holds true for $\varphi \in [\frac{\pi}{3}, \frac{\pi}{4}]$ what we promised.

Consequently (***) remains true if we lessen $|gh|$, so also if $0 < |gh| < 2\varphi(\delta) - \frac{1}{2}\pi$.
Still $|hm|$ increases as $|gh|$ decreases.
This ends the considerations of (3.3.3.c'') and thus of (3.3.3.c).

\vskip0.1cm
In (4) we act analogously as in (3).

\vskip0.1cm
The inequality in our theorem cannot be improved for every $\delta \in (\frac{1}{2}\pi, \pi)$.
This follows from the example of the regular triangle of sides of length $2\varphi(\delta)$, which is fulfilled when $|jf| = |jh| = 2\varphi(\delta)$ in Case 2.
\end{proof}

\vskip0.2cm
\section{Final remarks}

Formulating and proving an $n$-dimensional version of Theorem 2 seems to be very difficult.
Even formulating a conjecture does not seem to be easy.

From Theorem 2 we obtain the following corollary.

\vskip0.2cm 
\noindent
{\bf Corollary.} {\it For every convex body $C \subset S^2$ with $\diam(C) \in (\frac{\pi}{2}, \pi)$ we have} 
\vskip-0.2cm 
$${\rm diam}(E(C)) > \frac{2}{3} \cdot {\rm diam}(C).$$

\begin{proof}
By Theorem 2 it is sufficient to show that $\arccos \big(\frac{1}{4}(\cos \delta + \sqrt {\cos^2 \! \delta +8})\big) > \frac{1}{3}\delta$ for every $\delta \in (\frac{1}{2}\pi, \pi)$.
This is we have to show that 

$$\cos \delta + \sqrt {\cos^2 \! \delta +8} < 4 \cos \frac{1}{3}\delta \eqno (***\; *)$$  

\noindent
for every such a $\delta$.

Substituting $3\alpha = \delta$ into the well known formula $\cos 3\alpha = 4\cos^3\alpha - 3\cos\alpha$, we get $y= 4z^3 -3z$, where $y=\cos\delta$ and $z=\cos\frac{\delta}{3}$ (clearly, $z>0$).
Consequently, (****) is equivalent to $4z^3 -3z + \sqrt{(4z^3-3z)^2 +8} < 4z$, this is to $4z^4-5z^2+1 < 0$, which holds true for $\frac{1}{2} < z < 1$.
Hence (****) holds true if $\cos \frac{1}{6} \pi < \cos \frac{1}{6} \pi < \cos 0$, this is if $0< \delta < \pi$.
In particular this is true when $\frac{1}{2}\pi < \delta < \pi$.
\end{proof}

\vskip0.1cm
The inequality in Corollary cannot be improved.
In order to see this consider the limit of regular triangles whose heights tend to $\pi$ from below. 
The author expects that the following generalization of Corollary (where a convex body and its extreme points are defined as in $S^2$) holds true: {\it for every convex body $C \subset S^n$ of diameter $\delta  \in (\frac{\pi}{2}, \pi)$ we have} ${\rm diam} (E(C)) > \frac{2}{\pi} (\arcsin \sqrt \frac{n+1}{2n})\cdot \delta$. 
This value cannot be lessened for the family of the regular simplexes inscribed in the ball of radius $r$ as $r\to \pi/2$ from below.

\vskip0.2cm
{\bf Declarations}

{\bf Conflict of interests} The author has not disclosed any competing interests

{\bf Funding} The author declares no financial support for this work

{\bf Availability} This article has no associated data

\baselineskip 14pt

\vskip0.2cm
\noindent
Marek Lassak

\noindent
University of Science and Technology

\noindent
85-789 Bydgoszcz, Poland

\noindent
e-mail: lassak@pbs.edu.pl


\begin{thebibliography}{12}

\baselineskip=12pt

\bibitem{BS}
 B\"or\"oczky, K. J., Sagmeister, \'A., {\it Stability of the isodiametric problem on the sphere and in the hyperbolic space.} Adv. in Appl. Math. 145, 102480 (2023)

\bibitem {DGK} 
Danzer, L., Gr\"unbaum, B., Klee, V., Helly's theorem and its relatives, \textit{Proc. of Symp. in Pure Math.} vol. VII, Convexity, pp. 99--180 (1963) 

\bibitem{Ha} 
Hadwiger, H., Kleine Studie zur kombinatorischen Geometrie der Sph\"are.
{\it Nagoya Math. J.} 8, 45--48 (1955)

\bibitem{Han} 
Han, H., From spherical separation center set to the upper and lower bound theorems. arXiv:013797.v2

\bibitem{L1} 
Lassak M., Width of spherical convex bodies. {\it Aequ. Math.} 89, 555--567 (2015) 

\bibitem{L2} 
Lassak, M., Spherical Geometry -- A  Survey on Width and Thickness of Convex Bodies in {\it Surveys in Geometry.} I. Papadopoulos, Athanase (ed.), Cham: Springer, 7--47 (2022)

\bibitem{LasMus-AEQ}  
Lassak, M., Musielak M., Spherical bodies of constant width, {\it Aequ. Math.} 92, 627--640 (2018)

\bibitem{LasMus-Bull}  
Lassak, M., Musielak, M., Reduced spherical convex bodies, {\it Bulletin Polish Acad. Sci. Math.} 66, 87--97 (2018)

\bibitem{LasMus-FAS}  
Lassak, M., Musielak, M., Diameter of reduced spherical convex bodies, {\it Fasciculi Mathematici} 61, 103--108 (2018)

\bibitem{Sa}
Santalo, L. A., Convex regions in the $n$-dimensional spherical surfaces, {\it Ann. Math.} 47, 448--459 (1946)

\bibitem{Wh}
Whittlesay, M. A, {\it Spherical Geometry and Its Applications}, CRC Press, Taylor and Francis Group, 2020 

\bibitem{VB} Van Brummelen, G., {\it Heavenly Mathematics. The forgotten art of spherical trigonometry.} Princeton University Press, Princeton (2013)

\end{thebibliography}
\end{document}